\documentclass{amsart}
\usepackage{amsmath}
\usepackage{amsthm}

\theoremstyle{plain}
\newtheorem{lemma}{Lemma}[section]
\newtheorem{cor}[lemma]{Corollary}
\newtheorem{thm}[lemma]{Theorem}
\newtheorem{prop}[lemma]{Proposition}

\theoremstyle{definition}

\newtheorem{examples}[lemma]{Example}
\newtheorem{q}[lemma]{Question}

\newtheorem{remark}[lemma]{Remark}

\newcommand{\ann}{\operatorname{\ensuremath {\rm ann}}}

\newcommand{\xym}{\ensuremath \xymatrix@1}

\newcommand{\J}{\operatorname{\ensuremath{\it J}}}

\newcommand{\D}{\ensuremath {\mathcal{D}}}

\newcommand{\Aut}{\ensuremath {\rm Aut}}

\newcommand{\Q}{\ensuremath \mathbb{Q}}

\newcommand{\CC}{\ensuremath \mathbb{C}}

\newcommand{\Z}{\ensuremath \mathbb{Z}}

\newcommand{\M}{\ensuremath \mathbb{M}}

\newcommand{\U}{\operatorname{\ensuremath {\it U}}}
\newcommand{\rad}{\operatorname{\ensuremath {\rm rad}}}

\newcommand{\characteristic}{\operatorname{\ensuremath {\rm char}}}
\newcommand{\diag}{\operatorname{\ensuremath {\rm diag}}}
\newcommand{\gr}{\operatorname{\ensuremath {\rm gr}}}

\title{Morphic and principal-ideal group rings}
\author{Thomas J. Dorsey}
\address{Department of Mathematics, Vassar College, 124 Raymond Ave, Poughkeepsie, NY 12604} 
\email{thomasjdorsey@gmail.com}

\begin{document}
\bibliographystyle{amsplain}

\begin{abstract}
We observe that the class of left and right artinian left and right morphic rings agrees with the class of artinian principal ideal rings.   
For $R$ an artinian principal ideal ring and $G$ a group, we characterize when $RG$ is a principal ideal ring;
for finite groups $G$, this characterizes when $RG$ is a left and right morphic ring.   
This extends work of Passman, Sehgal and Fisher in the case when $R$ is a field, and work of Chen, Li, and Zhou on morphic group rings.  
\end{abstract}

\maketitle

\section{Introduction}

Throughout this article, the term artinian ring will refer to a left and right artinian ring, and the term principal ideal ring will refer to a 
ring all of whose one-sided ideals are principal.  
Consider the following question:
\begin{q} \label{mainq}
Given a ring $R$ and a group $G$, when is the group ring $RG$ a principal ideal ring?  
\end{q}

The classical group algebra case of this question, when $R$ is a (commutative) field, and $G$ is an arbitrary group, was 
answered by Sehgal and Fisher in \cite{sehgalfisher} in case $G$ is nilpotent, and completed by Passman in \cite[Theorem 4.1]{passman}.  
With only minor changes (detailed in the appendix, below), the proof given works for division rings as well.  

On another, perhaps seemingly unrelated, topic, 
in \cite{morph}, Nicholson and Sanchez-Campos investigated the ``morphic" rings (rings which satisfy the dual of the first isomorphism
theorem).  In \cite[Example 36]{morph}, an example is given 
of a group ring $RG$, for which $R$ is artinian and left and right morphic, $G$ is a finite group, but 
for which $RG$ is not a morphic ring.  Motivated by this example, in \cite{morphgrp}, J. Chen, Y. Li, and Y. Zhou investigated the question of when a group ring is morphic.  
In \cite[Section 2]{morphgrp}, they prove some general theorems about morphicity of group rings.  
For instance, $RG$ left morphic implies that $R$ is left morphic and $G$ is locally finite (\cite[Theorem 2.1]{morphgrp}); on the other hand, 
$RG$ is left morphic if $RH$ is left morphic for each finite subgroup $H$ of $G$ (\cite[Theorem 2.4]{morphgrp}).  
After these general theorems, they classify when $RG$ is left morphic in a few special cases: specifically, when $R$ is either semisimple or $\Z_n$ for some $n$, and $G$ 
is a finite abelian group.  
In addition, they complete the 
case when $G = D_n$ is a dihedral group and $R = \Z_{p^r}$ where $p$ is prime and $r \ge 1$.  The general problem of determination of when $RG$ is morphic is left open, even in the case when 
$R$ is a left and right artinian left and right morphic ring, and $G$ is a finite group.      
     
Nicholson and Sanchez-Campos also investigated the interplay between (left, right, or left and right) 
morphic rings and (left, right, or left and right) principal ideal rings in \cite{morphp}.  In particular, \cite[Corollary 16]{morphp} contains a structure 
theorem describing rings which are left artinian and left and right morphic, and it is shown that this class agrees with the class of rings 
which are left and right principally morphic.    
It seems to have been overlooked, however, that the structure theorem of \cite[Corollary 16]{morphp} (adding to the earlier \cite[Theorem 35]{morph}) 
contains a condition equivalent to the classical structure theorem describing the artinian principal ideal rings
found, for instance, in \cite[Section 15]{jacobsontheory} (stated with more modern terminology, for instance, in \cite[Corollary 2.2]{eisenbudgriffith}).
Namely, a ring is an artinian principal ideal ring if and only if it is a finite direct product of matrix rings over local artinian principal ideal rings.
In \cite[Corollary 16]{morphp}, it is shown that the left artinian left and right morphic rings are precisely the finite direct products of matrix 
rings over left and right ``special'' rings (in the terminology of \cite{morph}).  A left special ring is a local ring $R$ for which 
the Jacobson radical of $R$ is a left principal ideal, generated by a nilpotent element.  
It is easy to see that a ring is left special if and only 
if it is a local left artinian left principal ideal ring (using conditions (2) and (3) of \cite[Theorem 9]{morph} and the fact that a left artinian 
ring has a nilpotent Jacobson radical), and hence 
a left and right special ring is precisely a local artinian principal ideal ring.          

In view of this, there is another even more surprising equivalent condition that 
can be added to \cite[Corollary 16]{morphp}.  Namely, the class described there is, in fact, the artinian principal ideal rings, allowing us 
to add the first condition below.             
\begin{thm} $($cf. \cite[Corollary 16]{morphp}$)$ \label{princrings}  For any ring $R$, the following are equivalent:
\begin{enumerate}
	\item  $R$ is an artinian principal ideal ring. 
	\item  $R$ is left and right $P$-morphic.
	\item  $R$ is left artinian and left and right morphic.
	\item  $R$ is semiprimary and left and right morphic.
	\item  $R$ is left perfect and left and right morphic.
	\item  $R$ is a semiperfect, left and right morphic ring in which $J$ is nil and $S_r \subseteq^{ess} R_R$.
	\item  $R$ is a semiperfect, left and right morphic ring with ACC on principal left ideals in which $S_r \subseteq^{ess} R_R$.  
	\item  $R$ is a finite direct product of matrix rings over local artinian principal ideal rings $($i.e. left and right special rings$)$.  
\end{enumerate}
\end{thm}           
Consequently, whenever $R$ is a left and right artinian principal ideal ring and $G$ is a finite group, the 
group ring $RG$ is morphic if and only if $RG$ is a principal ideal ring.  
Thus, there is overlap in the study of morphic groups rings found in \cite{morphgrp} with the earlier study of principal ideal group rings found in 
\cite{sehgalfisher} and \cite[Section 4]{passman}.  
In particular, one of the special cases handled in 
\cite{morphgrp}, when $R$ is semisimple, and $G$ is a finite abelian group, is essentially already contained in the results of \cite[Section 4]{passman} and 
\cite{sehgalfisher}.  As we mentioned above, Passman, Seghal and Fisher only deal with the classical case with coefficients in a field, but 
their proofs essentially work in the case of a division ring, and easily imply a classification in the case of semisimple rings (see Theorem \ref{semisimple}, below).     

In this article,     
we will first answer Question \ref{mainq} in the case when $R$ is a local artinian principal ideal ring and $G$ is an arbitrary group.  
This extends \cite[Theorem 4.1]{passman}; simultaneously it includes as special cases observations made in \cite{morphgrp} on morphicity when the 
coefficient ring is $\Z_{p^n}$ (removing any hypothesis on the group $G$).      
We will then answer Question \ref{mainq} in 
the case when $R$ is an artinian principal ideal ring and $G$ is an arbitrary group.  
In particular, when $R$ satisfies the hypotheses of Theorem \ref{princrings} (e.g. is left artinian and left and right morphic), and $G$ is finite, 
we completely characterize when $RG$ is morphic.  In particular, our results contain each of the special cases dealt with in \cite[Section 3-4]{morphgrp}, and 
answer many of the question contained therein.      

We will, of course, rely heavily on \cite[Theorem 4.1]{passman} (for division rings), and this result will give us an extremely good start on our way.  
We will freely use the fact that all instances of ``field'' in \cite[Section 4]{passman} can be replaced by ``division ring'' (we detail this in the appendix, below).  
Also, Lemma \ref{lemma1}, below, is motivated by \cite[Theorem 2.8]{morphgrp}; aside from this motivation, we will not 
rely upon any of the results found in \cite{morphgrp}.   

Our ring-theoretic terminology will generally follow \cite{fc}.  In particular, for a ring $R$, we denote by $J(R)$ the Jacobson radical of $R$, and for 
a group ring $RG$, we denote by $\epsilon$ the augmentation map $\epsilon: RG \longrightarrow R$, whose kernel is the augmentation ideal $\Delta(RG)$.
For an element $x$ in a ring $R$, we denote by $\ann_{\ell}^R(x)$ and $\ann_r^R(x)$ the left and right annihilators of $x$ in $R$, respectively.  When the 
ring is clear from the context, we shall omit the superscript $R$.        
Also, a local ring is a (not necessarily noetherian) ring with a unique left (equiv. right)
ideal, which agrees with its Jacobson radical.  We shall also need some group-theoretic terminology.  
Specifically,  
if $\mathcal{A}$ and $\mathcal{B}$ are two classes of groups, we say that 
a group $G$ is $\mathcal{A}$-by-$\mathcal{B}$ if there exists $N \lhd G$ such that $N \in \mathcal{A}$ and $G/N \in \mathcal{B}$.  
Recall also that for finite groups $G$, we say that $G$ is a $p$-group if $|G|$ is a power of $p$, and we say that $G$ is a $p'$-group if $|G|$ is relatively prime to $p$.  
We will also allow ourselves the natural generalizations of this when $\pi$ is a finite set of primes.  In particular, if $\pi = \emptyset$, a finite $\pi'$-group is synonymous 
with a finite group, and the only finite $\pi$-group is the trivial group.  
We shall also freely use the fact that if $R$ is a local artinian principal ideal ring, then $\J(R) = cR = Rc$ for any $c \in \J(R) \setminus \J(R)^2$ 
(e.g. \cite[Corollary 10]{morph} or \cite[Theorem 38]{jacobsontheory}).   

\section{The local case}

As we mentioned in the introduction, Question \ref{mainq} has a complete characterization in the case when $R$ is a division ring.  Completing work of 
Fisher and Sehgal (\cite{sehgalfisher}) for nilpotent groups, in \cite[Section 4]{passman}, Passman showed that 
\begin{thm} \label{passman} \cite[Theorem 4.1]{passman}  Let $KG$ be the group ring of $G$ over the division ring $K$.  Then, the following are equivalent.  
\begin{enumerate}
	\item $KG$ is a right principal ideal ring.
	\item $KG$ is right Noetherian and the augmentation ideal $\Delta(KG)$ is principal as a right ideal.
	\item 
	\begin{itemize}
	\item		$\characteristic K = 0$:  $G$ is finite or finite-by-infinite cyclic.  
	\item 	$\characteristic K = p > 0$: $G$ is finite $p'$-by-cyclic $p$, or finite $p'$-by-infinite cyclic.  
	\end{itemize}
\end{enumerate}
\end{thm}  
As we mentioned above, the statement found in \cite[Theorem 4.1]{passman} requires that $K$ is a field, but with only minor changes (detailed in the appendix, below), its proof is valid when $K$ is a division ring 
as well.  For simplicity, we will refer to \cite[Theorem 4.1]{passman}, even in the case of division rings, as 
Passman's Theorem.

Before stating our extension of Passman's Theorem to local artinian principal ideal rings, 
we shall need a definition, which will require a bit of preliminary set up.  This work, and some of the 
work done when we discuss associated graded rings 
(in the beginning of Section \ref{assoc})  
is similar to that found in \cite[Chapter 2, Section 6]{jateg}, specifically, \cite[Chapter 2, Lemma 6.2]{jateg}.  We are working in a more restricted case 
when compared with that studied in \cite{jateg}, and a more elementary exposition is therefore possible.  In the interests of keeping the exposition elementary and relatively self-contained, we will 
deal explicitly with our special case, instead of extracting it from the results from \cite{jateg}.    

Suppose that $R$ is a local artinian 
principal ideal ring with $J^2=0$, and for which $J \ne 0$.  In this case, the only nontrivial ideal of $R$ is $J = cR = Rc$ (see, for instance, \cite[Theorem 9, ff.]{morph}).  
We associate to a ring with these properties a ring automorphism $\varphi$ of $R/J$ as follows.  Note that, for each $r \in R$, 
$cr = sc$ for some $s \in R$.  Since $\ann_{\ell}(c) = J$ (see \cite[Theorem 9]{morph}), it is clear that $s$ is determined uniquely as an element of $R/J$.    
Observe that $c1 = 1c$, and if $cr = sc$ and $cr' = s'c$, then $c(r+r') = (s+s')c$ and $c(rr') = scr' = (ss')c$.  
Thus, we have a well-defined ring homomorphism $\varphi: R \longrightarrow R/J$ defined by setting $\varphi(r) = s+J$, such that $cr = sc$.
Furthermore, 
$\sigma$ is surjective, since if $s \in R$, $sc \in Rc = cR$, so $sc = cr$ for some $r \in R$, and hence $\varphi(r) = s$.  
Observe that $cr = 0 = 0c$ if and only if $r \in \ann_r(c) = J$, so $\ker(\varphi) = J$.
We conclude that $\sigma$ induces a ring automorphism, which we will refer to as $\sigma$, of $R/J$.
Note that if $Rc' = c'R$, then $c' = uc$ for some $u \in \U(R)$.  Then, $c' r = uc r = us c = usu^{-1} uc = usu^{-1} c'$.  
Thus, the map $\sigma$ is only determined up to conjugation by a nonzero element of $R/J$.  For future reference, let us observe that 
the skew polynomial rings $(R/J)[t; \sigma]$ and $(R/J)[x; \rho_u \circ \sigma]$, where $\rho_u$ denotes conjugation by $u$, 
are isomorphic, by sending $t$ to $ux$.  
For our uses later, the possible conjugation of $\sigma$ will not be relevant (we will be dealing with skew polynomial rings, as above), 
so we will, in general, refer imprecisely to a single map $\sigma$.    

Given a local artinian principal ideal ring $R$, we associate to $R$ the (class of) $\sigma \in \Aut(R/J)$ corresponding to the 
above construction for the quotient ring $R/J^2$.  Now, if $G$ is a finite group $G$ with $|G| \cdot 1 \in \U(R/J)$, by 
Maschke's theorem, $(R/J)G$ is semisimple.  Let $\{e_1, \ldots, e_n\}$ be the set of centrally primitive idempotents of $(R/J)G$.  
The automorphism $\sigma$ extends to an automorphism of $(R/J)G$, acting on $G$ trivially, and 
must permute the set $\{e_1, \ldots, e_n\}$ of centrally primitive idempotents.  Note that this is reminiscent of Lemma 5 and Lemma 6 of \cite{sehgalfisher}.   
  
We shall say that a finite group $G$ with $|G| \cdot 1 \in \U(R/J)$ is $R$-admissible if $\sigma$ induces the identity permutation 
on the set of centrally primitive idempotents of $(R/J)G$.  Note that this does not depend on the choice of $\sigma$, since conjugation by 
a nonzero element of $R/J$ certainly must fix any central element of $(R/J)G$.    
The condition that $G$ is $R$-admissible is equivalent to saying that if $f_i$ is any lift of $e_i$ to $(R/J^2)G$, then 
$f_i c = c f_i$.  In particular, $R$-admissibility is actually the statement that the block decomposition of the artinian ring $(R/J)G$ lifts to a 
block decomposition of the artinian ring $RG$ (see \cite[Section 22]{fc}).  

In treating the group ring case specifically, we prefer to view $R$-admissibility as a property of the automorphism $\sigma$.
For instance, in the case when $k = R/J$ is an algebraically closed field, this is equivalent to the condition 
that $\sigma$ fixes $\chi(g)$ for each irreducible $k$-character $\chi$ of $G$, and each $g \in G$.  This is certainly ensured if $\sigma$ fixes all $|G|$-th roots of unity.  	 
Generalizing our results to other classes of rings may well be possible, however, using the lifting of block decompositions as one's starting point.    

Our main theorem for local artinian principal ideal rings is the following, which extends Passman's Theorem.  We do not know, however, if there is a 
valid analogue of condition (b) of Passman's Theorem.    
\begin{thm} \label{mainthm}
Suppose $R$ is a local artinian principal ideal ring and $G$ is a group.  Then, the following are equivalent:
\begin{enumerate}
	\item  $RG$ is a principal ideal ring 
	\item 
\begin{itemize}
	\item $\characteristic(R/J) = 0$: $G$ is a finite or finite-by-infinite cyclic. 
	If $R$ is not a division ring, then $G$ is an $R$-admissible finite group.  
	\item $\characteristic(R/J) = p>0$:  $G$ is finite $p'$-by-cyclic $p$, or a finite $p'$-by-infinite cyclic.  
	If $R$ is not a division ring, then $G$ is a finite $R$-admissible $p'$-group.  
\end{itemize}
\end{enumerate}
\end{thm}

Much of the forward implication follows immediately from Theorem \ref{passman}, since if $RG$ is a principal ideal ring, then $(R/J)G$ is as well, where $R/J$ is a division ring, so we may apply 
Theorem \ref{passman} to obtain information about the group $G$.  This gets us off to a very good start, however, there is still much to be done.      
In the forward implication, it remains only to show that, if $R$ is not a division ring, then $G$ is finite with $|G| \cdot 1 \in \U(R)$, and that 
$G$ is $R$-admissible.  Our next lemma will complete everything in the forward implication except for the $R$-admissibility.  
For the reverse implication, the entire case when $R$ is not a division ring 
remains.  We shall break the proof of Theorem \ref{mainthm} into a few steps, over the course of the next few sections.

Our first step is to prove two lemmas.  The first, in characteristic $p$, is motivated by 
\cite[Theorem 2.8]{morphgrp}.  
The argument in \cite[Theorem 2.8]{morphgrp} is specific to $\Z_{p^r}$ (possibly able to be extended to 
local artinian principal ideal rings for which $\J(R)$ has a central generator).  Our argument is completely different, obtaining a slightly weaker conclusion than 
\cite[Theorem 2.8]{morphgrp}, but for general local artinian principal ideal rings.    
When restricting to groups which admit surjections onto nontrivial $p$-groups whenever $p$ divides $|G|$ (e.g. finite 
nilpotent groups), we are able to obtain the same type of conclusion found in \cite[Theorem 2.8]{morphgrp}.      

\begin{lemma} \label{lemma1}  
Suppose $R$ is a local artinian principal ideal ring for which $\characteristic(R/J) = p$, and suppose $G$ is a finite $p$-group.  
If $RG$ is a principal ideal ring, then $R$ is a division ring.    
\end{lemma}
\begin{proof}
Suppose that $\characteristic(R/J(R)) = p> 0$, $G$ is a finite $p$-group, and that $RG$ is a principal ideal ring.    
By \cite{nichlocal}, $RG$ is an (artinian) local ring, which, by assumption is a principal ideal ring.    
It is easy to see that $\epsilon^{-1}(J(R))$ is the maximal left ideal of $RG$, since the left ideals of $RG$ form a chain (see \cite[Theorem 9]{morph}) and the only left ideals of $R$ are powers of $J$, so 
$I = \epsilon^{-1}(J(R))$ is a maximal left ideal of $RG$, and hence the unique maximal left ideal 
because $RG$ is local.  Also, it is apparent that $I^2 = \epsilon^{-1}(J(R)^2)$.          
In particular, if $R$ is not a division ring, then $\J(R) \setminus \J(R)^2$ is nonempty, and if $c \in \J(R) \setminus \J(R)^2$, the element $c \cdot 1$ is clearly an element of 
$I \setminus I^2$, hence it generates $I$ as a right ideal (e.g. \cite[Theorem 38]{jacobsontheory} or \cite[Claim 1, p. 395]{morph}).  But, since $G$ is nontrivial, we may find $1 \ne g \in G$, and the element $1 - g$ is an element of $I$.  
But $I = c(RG)$, so 
there must exist $x \in RG$ such that $(c \cdot 1) x = (1-h)$.  Comparing constant coefficients (noting that $c = c \cdot 1$ is a scalar), 
we find that $c$ is right invertible, which is clearly impossible, since $c \in \J(R)$.  
We conclude that $\J(R) = \J(R)^2$ so $\J(R) = 0$ (since $R$ is artinian), and hence $R$ is a division ring. 
\end{proof}
Our next lemma is in the same vein, for the group $\Z$, and applies to all characteristics.  
The previous lemma was stated only in the local case, since this is the only case we shall use, and since it is simpler to state due to the restriction on the characteristic.  
The following lemma will be just as easy to state without the condition that $R$ is local.  
The crux of the proof, however, is the local case, as in the last lemma.   
\begin{lemma}  \label{noz} Let $R$ be an artinian principal ideal ring.  If the ring $R \Z$, which is isomorphic to the Laurent polynomial ring $R \langle x \rangle$, 
is a principal ideal ring, then $R$ is semisimple.  
\end{lemma}
\begin{proof}   
Using the structure theorem for artinian principal ideal rings, we write $R \cong \prod_{i=1}^n \M_{k_i}(S_i)$, where $S_i$ is a local artinian principal ideal ring, and $k_i>1$.  
If $R$ is not semisimple, there is some $S_i$ for which $S_i$ is not a division ring.  Since $R$ is a principal ideal ring, so is its quotient $\M_{k_i}(S_i/J(S_i)^2)$.  
Thus, it suffices to consider the case when $R = \M_n(S)$, where $n \ge 1$ and $S$ is a local artinian principal ideal ring for which $\J(S)^2 = 0$.  

Let $T = R\langle x \rangle$, $\J(R) = Rc = cR$ and let $K = J(R) \langle x \rangle = c R \langle x \rangle = R \langle x \rangle c$, which is an ideal of $T$.
We shall write $\overline{T} = S/K \cong (R/J(R)) \langle x \rangle$ and for $t \in T$, we will denote by $\overline{t}$ the image of $t$ in $\overline{T}$.   
Since $R$ is artinian, $K \subseteq \J(T)$ by \cite[Proposition 9]{connell}.  
Consider the right ideal $I = (1+x)T + cT$.  We will show that $I$ is not principal.  
Note first that $\overline{I} = (1+x) \overline{T}$, which is a proper right ideal of $\overline{T}$.  We conclude that $I < T$.  
Suppose that $I = f T$.  We have $f = \sum_{i \in \Z} a_i x^i$.  Separate those coefficients which are in $\J(R)$ from those which are not, and write 
$f = f_0 + f_1$, where $f_1 \in \J(R)\langle x \rangle$, and each coefficient of $f_0$ is in $R \setminus \J$.  
Observe that $\overline{f} = \overline{f_0}$.    

Since $f T = I$, we have $fg = 1+x$ for some $g \in T$.  Note that $\overline{f_0} \overline{g_0} = 1 + x$, 
and that $1+x$ is not a zero divisor in $\overline{T}$, so $\ann_{\ell}(\overline{f_0}) = 0$.  
Also, $\overline{T} \cong (R/J(R)) \Z \cong \M_n(S/\J(S)) \Z$ is a left and right Noetherian ring, and is 
a prime ring by \cite[Theorem 10.20]{fc} and \cite[Connell's Theorem, p. 161]{fc}.  
We conclude that $\overline{T}$ is left and right nonsingular (e.g. \cite[Corollary 7.19]{lmr}) and hence we conclude that $\ann_r(\overline{f_0}) = 0$, by 
\cite[Lemma 10.4.9]{passmanbook}.

On the other hand, since $c \in I$, there must exist $h \in T$, written as $h = h_0 + h_1$, such that $fh = c$.  
Then, $c = fh = f_0 h_0 + f_1 h_0 + f_0 h_1$, since $\J(R)^2=0$.  Reducing modulo $K$, we see that 
$0 = \overline{c} = \overline{f_0} \overline{h_0}$.  Since $\ann_r(f_0) = 0$, we conclude that $\overline{h_0} = 0$.
By our choice of $h_0$ and $h_1$, we see that $h_0 = 0$.    
Now, $K = R \langle x \rangle c$, so we write $h_1 = g_1 c$ for some $g_1 \in R \langle x \rangle$.  
Thus, we have $c = f_0 g_1 c$, so $(1 - f_0 g_1) c = 0$.  Using \cite[Theorem 9]{morph}, we see that every coefficient of $1 - f_0 g_1$ is in $\J(R)$ 
(looking first at coefficients with respect to $x$, and then at the entries of the matrices).  We conclude that 
$f_0 g_1 \in 1 + \J(R) \langle x \rangle$.  In particular, since $\J(R) \langle x \rangle \subseteq \J(T),$ we conclude that $f_0 g_1$ is a unit.  
We conclude that $f_0$ is right invertible.  It is easy to see that $T = R \Z$ is noetherian, hence Dedekind-finite, so we conclude that $f_0$ is a unit.  
Therefore, $I = f f_0^{-1} T$, where $f f_0^{-1} = 1 + f_1 f_0^{-1}$.  Since each coefficient of $f_1 f_0^{-1}$ is in $\J(R)$, and $\J(R)^2 = 0$, we see that 
$f_1 f_0^{-1}$ is nilpotent.  We conclude that $f f_0^{-1} = 1 + f_1 f_0^{-1}$ is unipotent, hence a unit.  We conclude that $I = T$, a contradiction.     
     
\end{proof}

The next step is to reduce to the case when $J(R)^2=0$.  
This is motivated by \cite[Chapter 22]{fc}; the reduction to the case of square zero radical is a standard technique (see \cite[p. 332]{fc}) in studying artinian rings.  
\begin{lemma} \label{lift} Suppose $R$ is a local artinian ring.  Then, $R$ is a principal ideal ring 
if and only if $R/J^2$ is a principal ideal ring.    
\end{lemma}
\begin{proof}
We need only prove the reverse implication, since the class of principal ideal rings is 
closed under homomorphic images.  Assume $R/J^2$ is a principal ideal ring.    
Thus, $J/J^2 = c(R/J^2) = (R/J^2)c$ for some $c \in R/J^2$.
Lift $c$ to an element of $J \setminus J^2$ in $R$.    
Thus, $J^2 + cR = J$ and $J^2 + Rc = J$.  Since $R$ is artinian, $J$ is nilpotent.  
Thus, by \cite[Theorem 23.16]{fc}, we conclude that $Rc = J$ and $cR = J$.  By \cite[Theorem 9]{morph}, $R$ is a 
principal ideal ring.       
\end{proof}

Using the structure theorem for artinian principal ideal rings, we can easily remove the assumption that $R$ is local from Lemma \ref{lift}.
\begin{cor}  \label{lift2} Suppose $R$ is an artinian ring.  Then, $R$ is a principal ideal ring if
and only if $R/J^2$ is a principal ideal ring. 
\end{cor}
\begin{proof}
Only the reverse implication needs proof.  Assume $R/J^2$ is a principal ideal ring, so that 
$R/J^2 = \prod_{i=1}^n \M_{k_i}(S_i)$.  By \cite[Theorem 22.9]{fc}, we may lift the centrally primitive
idempotents $\{e_1, \ldots, e_n\}$ corresponding to the previous product to a full set $\{f_1, \ldots, f_n\}$ of centrally primitive idempotents 
of $R$.  Note that $f_i R f_i / \rad(f_i R f_i)  \cong  e_i R e_i / \rad(e_i R e_i) \cong \M_{k_i}(S_i/\rad(S_i))$.  
Since $f_i R f_i$ is artinian, we conclude by \cite[Theorem 23.10]{fc} that $f_i R f_i \cong \M_{l_i}(K_i)$ for some 
local ring $K_i$.  But then, $\M_{l_i}(K_i/J(K_i)^2) \cong \M_{k_i}(S_i)$.  
The uniqueness asserted in \cite[Theorem 23.10]{fc} implies that $l_i = k_i$ and $K_i/J(K_i)^2 \cong S_i$.  Thus, 
$K_i$ is a local artinian ring for which $K_i/J(K_i)^2$ is a principal ideal ring.  By Lemma \ref{lift}, 
we conclude that $K_i$ is a principal ideal ring.  We conclude that $R \cong \prod_{i=1}^n \M_{k_i}(K_i)$ is a 
principal ideal ring.    
\end{proof}

At this point, we obtain our desired reduction.  
\begin{cor} \label{reduce2}
If $R$ is a local artinian principal ideal ring and $G$ is a finite group with $|G| \cdot 1 \in \U(R)$, then $RG$ is a principal ideal ring if and only if $(R/J^2)G$
is a principal ideal ring.    
\end{cor}
\begin{proof}
By \cite[Proposition 9]{connell}, $\J(R)G \subseteq \J(RG)$ in case $R$ is artinian or $G$ is locally finite (both of those conditions are true in our situation).  
In this case, however, we obtain equality.  To see this, note that 
$$\frac{RG}{\J(R)G} \cong \left( \frac{R}{\J(R)} \right) G,$$ which is semisimple, since $R/J(R)$ is semisimple ($J$-semisimple and artinian) and $G$ is a finite group with 
$|G| \cdot 1 \in \U(R/J(R))$.  We conclude from \cite[Ex. 4.11]{fc}, that $\J(RG) \subseteq \J(R)G$, and hence $\J(RG) = \J(R)G$.  

Now, only the reverse implication needs proof, as usual.  We have $(R/J^2)G \cong RG/(\J(R)^2 G) \cong RG/(\J(RG)^2)$.  Since $RG$ is artinian, the result now follows 
from Corollary \ref{lift2}.    
\end{proof}

\section{Associated graded rings and the case $J^2=0$} \label{assoc}

In this section, we will handle completely the case when $R$ is a local artinian principal ideal ring with $\J(R)^2 = 0$, and $G$ is a finite group with $|G| \cdot 1 \in \U(R)$, and we 
shall use this to complete the proof of Theorem \ref{mainthm}.  
We will first look at the simpler case when $R$ is an associated graded ring with respect to its Jacobson radical.   
As we shall see  the prototype for this type of ring is a skew polynomial ring of the form $D[t;\varphi]$, where $\varphi$ is a ring automorphism of the division ring $D$
(cf. \cite[Chapter 2, Section 6]{jateg}).  

\begin{lemma} \label{source}
Suppose $D$ is a division ring, and $\varphi$ is a ring automorphism of $D$.  
Then, $D[t;\varphi]$ is a principal ideal ring.  
\end{lemma}
\begin{proof}
Apply \cite[Theorem 1.2.9]{robson}, 
noting that since $\varphi$ is an automorphism, $D[t;\varphi]$ can be viewed as both a right and left skew polynomial ring.  
\end{proof}

Note that $D[t;\varphi]/(t^2)$ is a local artinian principal ideal ring with radical $(t)$.  Now, suppose instead that we 
start with any local artinian principal ideal ring $R$ for which $J(R)^2 = 0$.    
We may form the associated graded ring of $R$ with respect to the ideal $J(R)$, which in this case is 
$(R/J(R)) \oplus J(R)$, since $J(R)^2 = 0$.  Fix a $c \in R$ such that $J(R) = cR = Rc$, and an associated ring automorphism $\sigma: R/J \longrightarrow R/J$.  
From the definition of $\sigma$, it is easy to see that $\gr_J R \cong (R/J)[t; \sigma]/(t^2)$, with $t$ corresponding to $c$ (recall that the choice of $\sigma$ does not 
affect the isomorphism type of $(R/J)[t;\sigma]$).    
Thus, eluded to earlier, rings of the form 
$R = D[t; \varphi]/(t^2)$ are the general form of local artinian principal ideal rings with $J(R)^2 = 0$ for which $\gr_J R \cong R$.        
Note that the study of group rings over such rings are much easier to study than general local artinian principal ideal rings with $J^2=0$, since 
$D[t;\varphi]/(t^2) G \cong DG[t;\varphi]/(t^2)$, where $\varphi$ is the automorphism of $DG$ obtained by extending the automorphism $\varphi$ linearly, acting trivially on $G$.  
Now, $DG$ is semisimple, but, the difficulty is that $\varphi$ need not respect the blocks (central idempotents) of $DG$.  It is clear that the blocks are preserved precisely when 
$G$ is $D[t;\varphi]$-admissible, and we shall see that this is precisely the case when $DG[t;\varphi]$ is a principal ideal ring.    

Therefore, we shall now work to characterize when $RG$ is a principal ideal ring, in the case that $R = \gr_J R$ is a local artinian principal ideal ring with $J(R)^2=0$.   
First, we shall need to describe the automorphisms of $\M_n(D)$, where $D$ is a division ring, so that we may characterize the skew polynomial rings with coefficients 
in the simple artinian ring $\M_n(D)$.    

\begin{lemma}  \label{conj} Let $D$ be a division ring, and $n>0$.  Let $\varphi$ be any ring automorphism of $\M_n(D)$.  
Then, $\M_n(D)[t; \varphi] \cong \M_n(D[x;\sigma])$ for some ring automorphism $\sigma$ of $D$; moreover, 
$\M_n(D)[t;\varphi] / (t^2) \cong \M_n(D[x;\sigma]/(x^2)).$      
\end{lemma}
\begin{proof}
Our first goal is to describe the automorphisms of $\M_n(D)$.  
Let $\{e_{ij} \}$ denote the usual matrix units of $\M_n(D)$, and set $f_{ii} = \varphi(e_{ii})$.  
Since $\{e_{11}, \ldots , e_{nn} \}$ is a collection of orthogonal local idempotents which sum to $1$, it is easy to see that the same is true for 
$\{f_{11}, \ldots, f_{nn} \}$.  By \cite[Exercise 21.17]{fc}, there is a unit $u \in \M_n(D)$ and a permutation $\pi \in S_n$ such that 
$f_{i,i} = u^{-1} e_{\pi(i),\pi(i)} u$.  There is a unit $v \in \M_n(D)$ (a permutation matrix) such that $v^{-1} e_{ii} v = e_{\pi(i),\pi(i)}$.  
Then, $f_{ii} = (vu)^{-1} e_{ii} (vu).$  

Now, let us consider the automorphism $\psi = \rho_{(vu)^{-1}} \circ \varphi$, where, throughout this proof, $\rho$ denotes conjugation.  By our choice of $u$ and $v$, $\psi(e_{i}) = e_{i}$.  
Consider $f_{ij} = \psi(e_{ij})$.  Note that $e_r f_{ij} e_s = \psi(e_r e_{ij} e_s).$  
If $i \ne r$ and $j \ne s$, then $e_r f_{ij} e_s = 0$.  Thus, if we express $f_{ij} = \sum_{ij} a_{ij} e_{ij}$, for $a_{ij} \in D$, 
we see that $f_{ij} = a_{ij} e_{ij}$.  Since $\psi(e_{i}) = e_i$, we see that $a_{ii}=1$.  Note that $a_{ij} a_{jk} = a_{ik}$ for all $i,j,k$.  In particular,
$a_{ij} a_{ji} = 1$ so each $a_{ij}$ is a unit.  Consider the diagonal matrix $w = \diag(a_{11},a_{21}, \ldots, a_{n1})$, which is an invertible matrix with 
inverse $\diag(a_{11}, a_{12}, \ldots, a_{1n})$.  Consider $\rho_w \circ \psi$.  
Note that $\rho_w \circ \psi (e_{ij}) = e_{ij}$ for all $i,j$.  

Now, let $d \in D$, and consider $d' = \rho_w \circ \psi(dI_n)$.  
Note that $e_{ii} d' e_{jj} = (\rho_w \circ \psi)( e_{ii} dI_n e_{jj} )$, which is zero if $i \ne j$.  Thus, $d'$ is a diagonal matrix.  
Consider the permutation matrix $x = e_{12} + e_{23} + \cdots + e_{n-1,n} + e_{n,1}$, whose inverse is $x^{-1} = e_{21} + e_{32} + \cdots + e_{n,n-1} + e_{1n}$.  
For any diagonal matrix $z$, conjugation by $x$ applies a cyclic shift on the entries of $z$.  In particular, $x^{-1} (dI) x = dI$, and hence 
$x^{-1} d' x = d'$.  It follows that $d'$ is a diagonal matrix of the form $cI_n$ for some $c \in D$.  It is easy to see that $(\rho_w \circ \psi)|_DI_n$ is a ring 
automorphism of $DI_n$.  We shall refer to this automorphism of $D \cong DI_n$ as $\sigma$, and we shall also use $\sigma$ to denote the automorphism of $\M_n(D)$ 
obtained by applying $\sigma$ componentwise.  

It is easy now to see that $\sigma^{-1} \circ \rho_w \circ \psi$ is the identity map on $\M_n(D)$.  
Indeed, $\sigma^{-1} \circ \rho_w \circ \psi$ fixes each $e_{ij}$ and fixes $DI_n$ elementwise, from 
which it follows that it fixes $\sum_{ij} a_{ij} e_{ij} = \sum_{ij} (a_{ij} I_n) e_{ij}.$  
We conclude that $\varphi = \rho_{uvw^{-1}} \circ \sigma.$

Let $z = (uvw^{-1})^{-1}$.    
Now, consider the map $g: \M_n(D)[t;\varphi] \longrightarrow \M_n(D[x;\sigma])$ defined by 
embedding $\M_n(D)$ in $\M_n(D[x;\psi])$ and sending $t$ to $zx$.
Note that $\M_n(D[x;\sigma]) \cong \M_n(D)[x;\sigma]$, where the first $\sigma$ refers to the ring automorphism of $D$, and the 
second refers to the ring automorphism of $\M_n(D)$ it induces componentwise.  
  
Note that $t A = \varphi(A) t$.  Note that $g(tA) = zxA$, and 
$g(\varphi(A) t) = \varphi(A) zx$, but $z^{-1} \varphi(A)z = \sigma(A)$, so 
$g(\varphi(A) t) = z \sigma(A) x = z x A = g(tA)$.  
Since $\M_n(D)[t;\varphi] = \M_n(D)[t]/\langle \{tA - \varphi(A) t: A \in \M_n(D)\} \rangle$, 
we conclude that $g$ is a well-defined ring homomorphism.  

Suppose that $p(t) = A_0 + A_1 t + \cdots + A_n t^n$, then 
$$\begin{aligned} g(p(t)) &= A_0 + A_1 zx + A_2 (zx)^2 + \cdots + A_n (zx)^n 
\\&= A_0 + A_1 z x + A_2 z \sigma(z) x^2 + \cdots + A_n z \sigma(z) \cdots \sigma^{n-1}(z) x^n. \end{aligned}$$ 
In particular, we see that $g$ is bijective, since $z$ is a unit.  We conclude that $\M_n(D)[t;\varphi] \cong \M_n(D[x;\sigma])$.      

Finally, note that the ideal $(t^2)$ is the set of polynomials with $A_0 = A_1 = 0$, and $(x^2)$ is the set of polynomials in 
$\M_n(D)[x;\sigma]$ of the form $B_2 x^2 + B_3 x^3 + \cdots + B_n x^n$.  It is clear that $g((t^2)) = (x^2)$, from which it follows that 
$\M_n(D)[t;\varphi]/(t^2) \cong \M_n(D[x;\sigma]/(x^2))$.       
\end{proof}

The following proposition is now essentially obvious.  
\begin{prop} \label{thesplitcase}
Suppose $R=\gr_J R$ is a local artinian principal ideal ring with $J^2=0$ and $G$ is a finite group with $|G| \cdot 1 \in \U(R)$.  
If $G$ is $R$-admissible, then $RG$ is a principal ideal ring.  
\end{prop}
\begin{proof}
We have $R \cong \overline{R}[t;\sigma]/(t^2)$ as above.  
By Maschke's Theorem with the Artin-Wedderburn Theorem, $\overline{R} G \cong \prod_{i=1}^k \M_{n_i}(D_i)$ for division rings $D_i$ and $n_i > 0$.  
The automorphism $\sigma$ of $\overline{R}$ extends to an automorphism of $\overline{R}G$, which, by assumption, fixes the centrally primitive idempotents 
(which correspond to the direct product decomposition above).  In particular, $\sigma$ acts as the direct product of automorphisms $\sigma_i$ of $\M_{n_i}(D_i)$.  
It is straightforward to see that 
$$
\begin{aligned} RG  & \cong (\overline{R}[t;\sigma]/(t^2)) G \cong (\overline{R}G)[t, \sigma]/(t^2) 
\cong \prod_{i=1}^k \left( \M_{n_i}(D_i) [t_i; \sigma_i] / (t_i)^2 \right)
\\ &\cong \prod_{i=1}^k \left( \M_{n_i}(D_i[x_i; \psi_i] / (x_i)^2 ) \right). \end{aligned}$$
By Lemma \ref{conj}, Lemma \ref{source}, and the structure theorem for artinian principal ideal rings, we conclude that $RG$ is a principal ideal ring.     
\end{proof}  

At this point, we know that if $G$ is $R$-admissible, then $\gr_{JG}(RG) = (\gr_J R)G$ is a principal ideal ring.       
What is not clear is the nature of the relationship between $RG$ and $\gr_{JG}(RG)$.  
Our goal is a theorem of the type sought in \cite[Chapter 2, Section 7]{jateg}.  We seek to conclude that a ring is a principal ideal ring, knowing 
that its associated graded ring (with respect to its Jacobson radical) is a principal ideal ring (note that $P(R) = J(R)$ when $R$ is artinian).  In general, this type of question is difficult, but 
we have imposed strong chain conditions which help us.  Moreover,   
the main trick that we need is that the $R$-admissibility of $G$ allows us to lift the centrally primitive idempotents of $(R/J)G$ to centrally primitive idempotents of $RG$.  

One important special case of this type of result (lifting through the associated graded ring) is found in \cite[Proposition 7.7]{jateg}; the following is a special case of that result.    
\begin{lemma}  \label{liftgr} Suppose that $R$ is a local artinian ring.  Then, the following are equivalent: 
\begin{enumerate}
	\item $R$ is a principal ideal ring, 
	\item $\gr_J R$ is a principal ideal ring.  
\end{enumerate}
\end{lemma}
\begin{proof}
Any (one-sided) artinian local ring is completely primary, and its prime radical agrees with its Jacobson 
radical (e.g. \cite[Theorem 10.30]{fc} and the fact that the Jacobson radical of an artinian ring is nilpotent).
The result is then a special case of \cite[Proposition 7.7]{jateg} applied on the left and the right.  

The special case when $J^2=0$ lends itself to a simpler proof, since $\gr_J R$ takes on a particularly simple form.  
In particular, the ideal $J$ in $R$ is also an ideal (i.e. $0 \oplus J$) of $\gr_J R$, and it is 
easy to check that the ideal $J$ is (left, resp. right) principal in $R$ if and only if $0 \oplus J$ is (left, resp. right) principal in $S$.  
\end{proof}

\begin{prop} \label{prop3}
If $R$ is a local artinian principal ideal ring with $J^2=0$, and $G$ is a finite group with $|G| \cdot 1 \in \U(R)$, then $RG$ is a principal ideal ring if and only if $G$ is $R$-admissible.  
\end{prop}
\begin{proof}
Let $S = RG$.  Since $G$ is finite and $|G| \cdot 1 \in \U(R)$, we see that $\J(RG) = \J(R)G$.    
First, lift the set $\{e_1, \ldots, e_n\}$ of centrally primitive idempotents of $RG$ to orthogonal primitive idempotents $\{f_1, \ldots, f_n\}$
(for instance, by \cite[Corollary 21.32]{fc}).    
Note that $f_1 + \cdots + f_n$ reduces to $e_1 + \cdots + e_n = 1$, modulo $\J(RG)$, so $f_1 + \cdots + f_n$ is an idempotent unit, hence equals $1$.   
If $i \ne j$, note that $f_i S f_j \subseteq \J(S)$.  To see this, note that reducing modulo $\J(RG) = \J(R)G$, we obtain
$e_i ( (R/J)G) e_j = 0$, since $e_i,e_j$ are orthogonal and central.  

For the forward implication, suppose $S$ is an artinian principal ideal ring.  By Theorem \ref{princrings}, 
$S$ is morphic.  By \cite[Corollary 19]{morph}, $f_i S f_j = 0$ whenever $i \ne j$.  We conclude that each $f_i$ is a central idempotent of $S$, since
$1 - f_i = \sum_{j \ne i} f_j$, so $f_i S (1-f_i) = 0 = (1-f_i) S f_i$, from which we conclude that $f_i$ is central (by \cite[Lemma 21.5]{fc}).      
In particular, $f_i c = c f_i$.  We conclude that $\sigma(e_i) = \sigma( \overline{f_i} ) = \overline{f_i} = e_i$.  We conclude that 
$G$ is $R$-admissible.    

For the converse, suppose that $G$ is $R$-admissible.  We claim that each $f_i$ is central in $S$.  First, note that 
$f_i c = c f_i$ for each $i$, by hypothesis, since $G$ is $R$-admissible.  Also, note that, even without $R$-admissibility, for any $r \in R$, 
$c(f_i r) = c(r f_i)$ and $(f_i r)c = (r f_i)c$, since $f_i r - r f_i$ is in $\J = \ann(c)$, since $\overline{f_i} = e_i$ is central in $\overline{S}$.     
Let $r \in R$, and let $i \ne j$.  Note that $f_i r f_j \in f_i S f_j \subseteq \J(S)$.  
We may write $f_i r f_j = c h = h' c$, for some $h,h' \in S$, since $\J(S) = \J(R)G = c(RG) = (RG)c.$
Note that $f_i r f_j = f_i(f_i r f_j) = f_i ch = c f_i h = c h f_i = f_i r f_j f_i = 0$.  
We conclude that $f_i S f_j = 0$ if $i \ne j$, and, as before, that each $f_i$ is central.  
We conclude that the $f_i$ are a complete set of 
centrally primitive idempotents of $S$.
The ring $f_i S f_i$ is an artinian ring which is a simple artinian ring modulo its radical.  By 
\cite[Theorem 23.10]{fc}, $f_i S f_i \cong \M_{k_i}(S_i)$ for some $k_i > 0$ and some local ring $S_i$.  
At this point, we know that $RG \cong \prod_{i=1}^n \M_{k_i}(S_i)$.   
Next, we will show that each $S_i$ is a principal ideal ring.   

Note that $(\gr_J R)G \cong \gr_J (RG) \cong \prod_{i=1}^n \M_{k_i}(\gr_J S_i).$  
By hypothesis, $G$ is $R$-admissible, so it is also $\gr_J R$ admissible (the induced automorphism of $R/J$ is the same); by Proposition \ref{thesplitcase}
$(\gr_J R)G$ is an artinian principal ideal ring.  
Since the class of artinian principal ideal rings is Morita invariant (e.g. \cite[Corollary 17]{morphp}), 
we conclude that $\gr_J S_i$ is a local artinian principal ideal ring.  
By Lemma \ref{liftgr}, we conclude that $S_i$ is a local artinian principal ideal ring.  
It follows that $RG \cong \prod_{i=1}^n \M_{k_i}(S_i)$ is a principal ideal ring.  
\end{proof}
\begin{remark}  
In light of the Artin-Wedderburn Theorem and the structure 
theorem for artinian principal ideal rings, we view Proposition \ref{prop3} as an analogue of Maschke's Theorem (in the form of \cite[Theorem 6.1]{fc}).
\end{remark}      

We can now put all of this together to prove Theorem \ref{mainthm}.  
\begin{proof}(of Theorem \ref{mainthm})
For the implication $(1) \Longrightarrow (2)$, 
suppose $RG$ is a principal ideal ring.  Then, its quotient $(R/J)G$ is a principal ideal ring, to which we may apply Theorem \ref{passman}.    
If $\characteristic(R/J) = 0$, we conclude that $G$ is finite or finite-by-infinite cyclic.  In the latter case, $RG$ surjects onto 
$R\Z$, which must therefore be a principal ideal ring.  By Lemma \ref{noz}, we conclude that $G$ must be finite if $R$ is not a division ring.  

If $\characteristic(R/J) = p>0$, we conclude from Theorem \ref{passman} that $G$ is finite $p'$-by-cyclic $p$ or finite $p'$-by-infinite cyclic.  
In either case, if the cyclic group in question is nontrivial, $RG$ surjects onto $RH$ for some nontrivial cyclic $p$-group $H$, and $RH$ is a principal ideal ring.
By Lemma \ref{lemma1}, we conclude that, if $R$ is not a division ring, then the cyclic group must be trivial, and hence $G$ must be a $p'$-group.      

We have shown that if $R$ is not a division ring, then $G$ is a finite group with $|G| \cdot 1 \in \U(R)$.  Since $RG$ is a principal ideal ring, its homomorphic image
$(R/J^2)G$ is a principal ideal ring, and by Proposition \ref{prop3}, $G$ is $R$-admissible.   

For the implication $(2) \Longrightarrow (1)$, Theorem \ref{passman} handles the case
when $R$ is a division ring.  In the remaining case, $R$ is not a division ring, $G$ is an 
$R$-admissible finite group with $|G| \cdot 1 \in \U(R)$.  By Proposition \ref{prop3}, $(R/J^2)G$ is a principal ideal ring, and by  
Corollary \ref{reduce2}, we conclude that $RG$ is a principal ideal ring.      
\end{proof}

\section{General artinian principal ideal rings}

Given the structure theorem for artinian principal ideal rings, we are now in position to easily study $RG$ when $R$ is an artinian principal ideal ring.  
We shall first need a few completely elementary group theoretic lemmas.  
\begin{lemma} \label{silly} Suppose that $G$ is a finite group, and that $H,K \lhd G$ such that $G/H$ and $G/K$ are cyclic groups with relatively 
prime order.  Then, $G/(H \cap K)$ is a cyclic group of order $|G/H| \cdot |G/K|$.

Suppose that $k \ge 1$ and $H_1, \ldots, H_k$ are normal subgroups of a finite group $G$, such that $G/H_i$ is cyclic for each $i$, and that
$|G/H_i|$ is relatively prime to $|G/H_j|$ if $i \ne j$.  Then, $G/(H_1 \cap \cdots \cap H_k)$ is a cyclic group of order $|G/H_1| \cdot |G/H_2| \cdots |G/H_k|$.  

\end{lemma}
\begin{proof}
Write $|G| = rst$, where $|G/H| = r$, $|G/K| = s$.  Since $H,K \lhd G$, then $HK \lhd G$, and contains $H$ and $K$.  
In particular, its order is divisible by $|H| = st$ and $|K| = rt$, so $|HK| = rst$, since $(r,s) = 1$.  We conclude that $HK = G$.  
The groups $H/(H \cap K)$ and $K/ (H \cap K)$ intersect trivially and have product equal to $G/(H \cap K)$, so $G/(H \cap K) \cong H/(H \cap K) \times G/ (H \cap K)$.  
But, $H/(H \cap K) \cong HK/K \cong G/K$ and $K/(H \cap K) \cong HK/H \cong G/H$ 
are cyclic groups with relatively prime orders $r$ and $s$ respectively, so we conclude that $G/(H \cap K)$ is cyclic.  

For the next statement, we induct, using the first statement.  In case $n=1$ the result is trivial, and we have already done the $n=2$ case.  
Working by induction, $G/(H_1 \cap H_2 \cap \dots \cap H_{k-1})$ is cyclic with order $|G/H_1| \cdots |G/H_{k-1}|$.  
Applying the $n=2$ case to $(H_1 \cap H_2 \cap \cdots \cap H_{k-1})$ and $H_k$, we find that 
$G/(H_1 \cap H_2 \cap \cdots \cap H_k)$ is a cyclic group whose order is $|G/H_1| \cdots |G/H_k|$.     

\end{proof}

\begin{cor}  \label{ptopi} Let $\pi$ be a nonempty finite set of primes, and let $G$ be a group.  
Then, $G$ is finite $\pi'$-by-cyclic $\pi$ if and only if $G$ is finite $p'$-by-cyclic $p$ for each $p \in \pi$. 
\end{cor}
\begin{proof}
We are given that for each $p \in \pi$, there exists a normal $p'$-subgroup $H_p \lhd G$ such that $|G/H_p|$ is a cyclic $p$-group.  
Applying Lemma \ref{silly}, we find that the group $G/(\bigcap_{p \in \pi} H_p)$ is cyclic and its order is $\prod_{p \in \pi} |G/H_p|$.  In particular, 
$G/(\bigcap_{p \in \pi} H_p)$ is a cyclic $\pi$-group.  On the other hand, $\bigcap_{p \in \pi} H_p$ is a subgroup of each $H_p$, which is a $p'$-group.  
We conclude that $\bigcap_{p \in \pi} H_p$ is a finite normal $\pi'$-group, and we conclude that $G$ is finite $\pi'$-by-cyclic $p$.     
\end{proof}
Similarly, we have the following lemma in the infinite case.  
\begin{lemma}  \label{infcyc} Let $\pi$ be a nonempty finite set of primes, and let $G$ be a group.  Then, $G$ is finite $\pi'$-by-infinite cyclic if and only if 
$G$ is finite $p'$-by-infinite cyclic for each $p \in \pi$.  
\end{lemma}
\begin{proof}
Suppose that $G$ is finite $\pi'$-by-infinite cyclic, so there is a finite $\pi'$-group $H \lhd G$ such that $G/H$ is infinite cyclic.  For each $p \in \pi$, 
$H$ is a $p'$-group.  We conclude that $G$ is finite $p$-by-infinite cyclic.  

On the other hand, suppose that $G$ is finite $p'$-by-infinite cyclic for each $p \in \pi$.  
We claim that $G$ is finite $\pi'$-by-infinite cyclic.  
We induct on the size of $\pi$, the result being trivial if $|\pi| = 1$.  
Thus, suppose that $\pi = \pi_1 \cup \{p\}$, where $|\pi_1| < |\pi|$.  
By the inductive hypothesis, $G$ is finite $\pi_1'$-by-infinite cyclic.  Thus, we have a finite normal $\pi_1'$-subgroup $H_1 \lhd G$ such that 
$G/H_1$ is infinite cyclic.  We also have a finite normal $p'$-subgroup $H \lhd G$ such that $G/H$ is infinite cyclic.  
Note that the subgroup $H_1 H$ is a finite normal subgroup, and its image in $G/H$ and $G/H_1$ is thus trivial, since an infinite cyclic group has no nontrivial finite subgroups.  
We conclude that $H = H_1$.  In particular, $H = H_1$ is a normal $\pi$-subgroup for which $G/H$ is infinite cyclic.  We conclude that $G$ is finite $\pi'$-by-infinite cyclic.  
\end{proof}

We reformulate Passman's Theorem for semisimple rings as follows; note that the statement is somewhat simpler than that of Passman's Theorem, 
not distinguishing the characteristic.      
\begin{thm}  \label{semisimple}
Let $R$ be a semisimple ring, and $G$ a finite group.  Let $\pi$ be the set of primes which are not invertible in $R$.  Then, the following are equivalent:
\begin{enumerate}
	\item $RG$ is a principal ideal ring.
	\item $G$ is finite $\pi'$-by-cyclic $\pi$, or finite $\pi'$-by-infinite cyclic.    
\end{enumerate}
\end{thm}
The statement for artinian principal ideal rings is similarly the following.    
\begin{thm} \label{mainart}
Suppose $R$ is an artinian principal ideal ring, and that $G$ is a finite group.  
Write $R \cong \prod_{i=1}^n \M_{k_i}(S_i)$ where each $S_i$ is a local artinian principal ideal ring. 
Let $\pi$ denote the set of 
primes which are not invertible in $R$.    
Then, the following are equivalent:
\begin{enumerate}
	\item  $RG$ is a principal ideal ring.   
	\item  $G$ is finite $\pi'$-by-cyclic $\pi$ or finite $\pi'$-by-infinite cyclic.  
	If $R$ is not semisimple, then $G$ is finite, and 
	for each $i \in \{1, \ldots, n\}$ for which $S_i$ is not a division ring,  
	$|G| \cdot 1 \in \U(S_i)$ and $G$ is $S_i$-admissible.      
\end{enumerate}
\end{thm}

We shall prove Theorems \ref{semisimple} and \ref{mainart} together.  In case $G$ is finite, these theorems 
are simply repackagings of their analogues Theorem \ref{passman} and Theorem \ref{mainthm}, using the relevant structure theorems and the 
fact that the class of {\em artinian} principal ideal rings is Morita invariant (passing to and from matrix rings).  We shall need to do some work, 
however, even in the semisimple case, to deal with arbitrary (not necessarily finite) groups.  The arguments in that case model essentially those found in \cite[Section 4]{passman}, but need to be adapted 
slightly.  In particular, we do not know whether, even for semisimple rings, an analogue of condition $(b)$ of Passman's Theorem holds.         

\begin{proof}(Proof of Theorem \ref{semisimple} and \ref{mainart})
The reverse implication is straightforward.  
Indeed, either $S_i$ is a division ring and $G$ is finite $\pi'$-by-cyclic $\pi$ or finite $\pi'$-by-infinite cyclic; or else, 
$|G|$ is finite with $|G| \cdot 1 \in \U(S_i)$, and $G$ is $S_i$-admissible.  By Theorem \ref{mainthm} and Theorem \ref{passman}, 
$S_i G$ is a principal ideal ring.  By \cite[Theorem 40]{jacobsontheory}, $\M_{k_i}(S_i G)$ is a principal ideal ring, and by 
\cite[Lemma on p.70]{jacobsontheory}, $R \cong \prod_{i=1}^n \M_{k_i}(S_i)$ is a principal ideal ring.        

The forward implication requires a slight amount of work.  
As we mentioned before starting the proof, in case $G$ is finite, 
this work evaporates, since the class of artinian principal ideal rings is Morita invariant; in particular, 
if $RG$ is an artinian principal ideal ring, so is $S_i G$ for each $i$, from which we can easily complete the argument.  
For principal ideal rings we do not know, in general, whether the (full) Peirce corner rings of a 
principal ideal ring need to be principal ideal rings.  We shall sidestep this problem, however.    

First, let us deal with the semisimple case.  
We will suppose first that $R \cong \M_n(K)$ is simple artinian, where $n \ge 1$ and $K$ is a division ring.  
Essentially, we will argue as in the proof of the implication $(b) \Longrightarrow (c)$ in 
Passman's Theorem (to obtain information about $G$), however, we will need to adapt those arguments slightly to our situation.  We will, however, 
only prove the implication $(a) \Longrightarrow (c)$ in this situation, which allows us more flexibility.       

Proceeding as in the implication $(b) \Longrightarrow (c)$ (though we are doing the analogue of the implication $(a) \Longrightarrow (c)$) of the proof of 
Passman's Theorem, we conclude that $RG$ is noetherian, and hence all subgroups of $G$ are finitely generated, and in particular, we have $\Delta^{+}(G)$ finite.  
As in Passman's Theorem, setting $\overline{G} = G/\Delta^+(G)$, we see that $K\overline{G}$ is a prime ring (e.g. \cite[Connell's Theorem, p. 161]{fc}), and 
hence $RG \cong \M_n(K)\overline{G} \cong \M_n(K \overline{G})$ is a prime ring (by \cite[Theorem 10.20]{fc}).  

At this point, we seek to apply \cite[Lemma 4.4]{passman}, which doesn't apply directly to our situation.  
Fortunately, we are only attempting to prove an analogue of the implication $(a) \Longrightarrow (c)$ as opposed 
to the more restrictive implication $(b) \Longrightarrow (c)$.  We will show that the conclusion of \cite[Lemma 4.4]{passman} is valid
if we assume $RG$ is a principal ideal ring, instead of only assuming that $\Delta(RG)$ is a principal right ideal and $RG$ is noetherian.  

Indeed, using the argument found in \cite[Lemma 4.4]{passman}, with $R$ simple artinian\footnote{In the appendix, below, we observe that the basic properties of the 
dimension subgroups needed apply in this situation, since $\Q$ or $\Z_p$ embeds in $R$.}, we conclude as before that when $\characteristic(K) = 0$\footnote{Note that in this case 
$\Q$ embeds in $R$, so we obtain the usual basic properties of the dimension subgroups; see the appendix, below.}, $G/G'$ is infinite cyclic.  
Similarly, if $\characteristic(K) = p>0$ we conclude that $|G/H|$ is infinite, where $H = \bigcap_{n=1}^{\infty} \D_n(RG)$.  
As before, we conclude that $G_1 = G/\D_n(RG)$ is a finite $p$-group.  Now, $R G_1 \cong \M_n(KG_1)$ is an {\em artinian} principal ideal ring, and, by Morita invariance, 
we conclude that $KG_1$ is a principal ideal ring.  Now, we are in position to apply \cite[Lemma 4.3]{passman} (for division rings, see the appendix), 
to conclude that $G/\D_n(RG)$ is cyclic, and hence $G' \subseteq H$, so 
$G/G'$ is infinite.  The rest of the proof of \cite[Lemma 4.4]{passman} carries through routinely.  
Indeed, $G/G'$ is an infinite finitely generated abelian group, so there is a normal subgroup
$W$ of $G$ for which $G/W$ is infinite cyclic.  If we set $B = \Delta(RW) RG$, then $B$ is a prime ideal of $RG$ since $RG/B \cong R(G/W) \cong \M_n(K(G/W))$ is a prime ring by 
\cite[Connell's Theorem, p.161]{fc} and \cite[Theorem 10.20]{fc}.  Applying \cite[Lemma 4.2(ii)]{passman}, we conclude that $B=0$, and hence $W = 1$.

Returning to the main proof, 
we therefore conclude that $\overline{G}$ is infinite cyclic or else $\overline{G} = 1$.  Thus, $G$ is finite or finite-by-infinite cyclic, and if $\characteristic(K) = 0$, we 
are done.  Now, suppose $\characteristic(K) = p>0$.  If $G$ is finite, then, since $RG \cong \M_n(KG)$ is an artinian principal ideal ring, we conclude that $KG$ is a principal ideal ring, 
and hence we conclude from Theorem \ref{passman} that $G$ is finite $p'$-by-cyclic $p$.  Finally, we argue as in the last paragraph of the implication $(b) \Longrightarrow (c)$ of 
the proof of Passman's Theorem, and, instead of applying \cite[Lemma 4.3]{passman} to $R \tilde{G}$, we note first that $R \tilde{G} \cong \M_n(K \tilde{G})$ 
is an artinian principal ideal ring, so $K \tilde{G}$ is a principal ideal ring, to which we may apply \cite[Lemma 4.3]{passman}, and we conclude as in the original proof, that 
either $G$ is finite $p'$-by-cyclic $p$, or else $G$ is finite $p'$-by-infinite cyclic.  

Putting all of the information we have together, suppose now that $R$ is semisimple, so that $R \cong \prod_{i=1}^n \M_{k_i}(K_i)$, and that 
$G$ is a group, for which $RG$ is a principal ideal ring.  Looking at quotients, we find that $\M_{k_i}(K_iG)$ is a principal ideal ring.  
Let $\pi$ be the set of primes which are not invertible in $R$; equivalently, $p \in \pi$ if and only if $p = \characteristic(K_i)>0$ for some $i$.      
First suppose that $G$ is finite.  
If $\pi$ is empty, then considering any $i$, we find that $G$ is finite; equivalently, $G$ is finite $\pi'$-by-cyclic $\pi$.  
If $\pi$ is nonempty, then, we find that $G$ is finite $p'$-by-cyclic $p$ for each $p \in \pi$ (considering any $i$ for which $K_i$ has characteristic $p$). 
By Lemma \ref{ptopi}, we conclude that $G$ is finite $\pi'$-by-cyclic $\pi$.  
We conclude in each case that if $G$ is finite, then $G$ is finite $\pi'$-by-cyclic $\pi$.  

Next, suppose that $G$ is infinite.  If $\pi = \emptyset$, then each $K_i$ has characteristic $0$; we conclude that $G$ is finite-by-infinite cyclic (we conclude 
this for each $i$); equivalently, $G$ is finite $\pi'$-by-infinite cyclic.  
If $\pi$ is nonempty, then, we conclude for each $p \in \pi$ that $G$ is finite $p'$-by-infinite cyclic
(considering any $i$ for which $K_i$ has characteristic $p$).  
By Lemma \ref{infcyc}, we conclude that $G$ is finite $\pi'$-by-infinite cyclic.  We therefore conclude that, in any event, $G$ is finite $\pi'$-by-infinite cyclic, and the forward implication has been proved when $R$ is semisimple.

With the semisimple case completed, we will now tackle the general case.  
Suppose that $RG$ is a principal ideal ring, where $R$ is an artinian principal ideal ring.  
Thus, $(R/J)G$ is a principal ideal ring, but $R/J$ is semisimple (since it is $J$-semisimple and artinian) 
and applying the semisimple case, we find that $G$ is finite $\pi'$-by-cyclic $\pi$ or finite $\pi'$-by-infinite cyclic.
In the latter case, $R \Z$ is a quotient of $RG$, so $R \Z$ is a principal ideal ring.  By Lemma \ref{noz}, we conclude that $R$ is semisimple.  In particular, if 
any $S_i$ is not a division ring, then $G$ is finite $\pi'$-by-cyclic $\pi$.    
Thus, suppose that $S_i$ is not a division ring.  Clearly, since $RG$ is a principal ideal ring, 
its quotient $\M_{n_i}(S_iG)$ is a principal ideal ring as well.  
Since $G$ is finite $\pi'$-by-cyclic $\pi$, $G$ is finite, so $\M_{n_i}(S_iG)$ is an artinian principal ideal ring, so $S_i G$ is a principal ideal ring. 
By Theorem \ref{mainthm}, we conclude that $|G| \cdot 1 \in \U(S_i)$ and $G$ is $S_i$-admissible, which completes the forward implication.  

For the implication $(1) \Longrightarrow (2)$, Theorem \ref{mainthm} implies that 
$G$ is a finite $p'$-by-cyclic $p$ group for each $p \in \pi$; if $S_i$ is not a division ring, then $|G| \cdot 1 \in \U(R)$ and $G$ is $S_i$-admissible.
By Corollary \ref{ptopi}, $G$ is a finite $\pi'$-by-cyclic $\pi$ group, and the proof is complete.  
\end{proof}

\begin{remark}
For $R$ a semisimple ring and $G$ a finite abelian group, \cite[Theorem 3.7]{morphgrp} characterizes when $RG$ is (strongly) left morphic, which, by Theorem \ref{princrings}
is equivalent to $RG$ is a principal ideal ring.  The second condition found there, that 
that for each $p \in \pi$, each Sylow $p$-subgroup of $G$ is cyclic, is equivalent to condition (2) above in the case when $G$ is a finite nilpotent group.

Also, Theorem \ref{mainart}, in the case when $R$ is a commutative artinian principal ideal ring (for which any group is $R$-admissible), reduces to the statement that 
$RG$ is a principal ideal ring if and only if $G$ is finite $\pi'$-by-cyclic $\pi$ and, for any $p \in \pi$, if $p \in \J(S_i)$, then $S_i$ is a division ring.  
The characterization of when $\Z_n G$ is morphic (i.e. a principal ideal ring) appears in \cite[Theorem 3.15]{morphgrp}, but its equivalence to this condition 
is somewhat obscured, since the statement and proof find the number theoretic condition that $p^2$ does not divide $n$, 
which happens to be equivalent to the aforementioned ring-theoretic condition for the ring $\Z_n$.    
\end{remark}

\section{Examples}

Our work settles a number of questions raised in \cite{morphgrp}.
In particular, Theorem \ref{mainart} answers in the affirmative \cite[Conjecture 4.14]{morphgrp} and \cite[Question 4.15]{morphgrp}.   
Next, let us resolve \cite[Question 2.6]{morphgrp} in the negative.  
First, we will need the following useful example, due to the author and A. Diesl.    

\begin{examples} \label{alex}
Let $R=\CC[t;\sigma]/(t^2)$, where $\sigma$ is complex conjugation, and $G = C_3$.  
Note that $G$ is not $R$-admissible, since the central idempotents of $(R/J)G \cong \CC C_3$ are $\frac{1 + \alpha g + \alpha^2 g^2}{3}$, where 
$\alpha^3 = 1$, but $\sigma$ does not fix the cube roots of unity.  
By Theorem \ref{mainthm}, $RG$ is not a principal ideal ring (equivalently, it is not morphic, since $RG$ is artinian).  
\end{examples}

Now, let us use Example \ref{alex} to answer \cite[Question 2.6]{morphgrp} in the negative.  
Consider the ring $R$ from Example \ref{alex}, let $H = C_3$, viewed as a subgroup of $G = S_3$.  
We have seen that $RC_3$ is not a principal ideal ring.  
The ring $RS_3$ is, however, a principal ideal ring, since 
$S_3$ is $R$-admissible (since each entry of the character table of $S_3$ is in $\Z$, the coefficients of the centrally primitive 
idempotents of $\CC S_3$ are all in $\Q$, and hence are fixed by $\sigma$).  

We note, in passing, that if $R$ is a local artinian principal ideal ring for which $\J(R)$ has a central generator, then every finite 
group $G$, for which $|G| \cdot 1 \in \U(R)$, is $R$-admissible, since the automorphism of $R/J$ is the identity map.  In particular, 
if $G$ is finite, $R$ is such a ring, then if $RG$ is morphic, the same is true for any subgroup $H$ of $G$ (since $|H|$ divides $|G|$).   

It should also be noted that the likely motivation for the authors of \cite{morphgrp} to ask \cite[Question 2.6]{morphgrp} lies in the 
statement and proof of \cite[Theorem 2.4]{morphgrp}.  The full strength of the hypotheses of \cite[Theorem 2.4]{morphgrp} are not 
needed in the proof and can be weakened.    
Namely, if $G$ is a locally finite group 
with the property that every element $x \in RG$ is left morphic as an element of $RH$ for some finite subgroup $H$ of $G$, then $RG$ is 
left morphic (instead of letting $H$ be the subgroup generated by the support, simply take $H$ to be the finite subgroup for which 
$x$ is left morphic in $RH$; the rest of the proof is unchanged).        
In fact, this gives a local condition for morphicity in the group ring $RG$.    
\begin{prop}  Let $G$ be a locally finite group.  Then, $RG$ is left morphic if and only if for each $x \in RG$, there is a finite subgroup $H$ of $G$ 
such that $x$ is left morphic in $RH$.     
\end{prop} 
\begin{proof}
The reverse implication is proved as in the proof of \cite[Theorem 2.4]{morphgrp} (replacing the subgroup generated by the support of the element by 
the larger subgroup guaranteed by hypothesis).  For the converse, if $x \in RG$ is left morphic, then there is some $y$ such that $\ann_{\ell}^{RG}(x) = RG y$ and 
$\ann_{\ell}^{RG}(y) = RG x$.  Consider the subgroup $H$ generated by the supports of $x$ and $y$; $H$ is finite, since $G$ is locally finite.  
Clearly, $RH x \subseteq \ann_{\ell}^{RH}(y)$ and $RH y \subseteq \ann_{\ell}^{RH}(x)$.  
Conversely, if $z \in \ann_{\ell}^{RH}(y)$, then $z \in \ann_{\ell}^{RG}(y) = RG x$.  
We write $z = w x$, for some $w \in RG$, and write $w = \sum g_i b_i$, where $\{g_i\}$ is a left transversal for $H$ in $G$, and $b_i \in RH$.  
We have $z = \sum g_i (b_i x)$.  Comparing coefficients of elements of $H$, we see that $z = g_0 b_0 x$, where $g_0 \in H$.
Since $g_0 b_0 \in RH$, we conclude that $z \in RH x$, and hence $RH x \supseteq \ann_{\ell}^{RH}(y)$.  It follows that $RH x = \ann_{\ell}^{RH}(y)$.  Similarly, 
$RH y = \ann_{\ell}^{RH}(x)$, and we conclude that $x$ is left morphic in $H$.      
\end{proof}

An interesting question, however, is 
whether $RG$ left morphic implies that for each $x \in RG$, there is a finite subgroup $H$ of $G$ for which 
$x \in RH$ and $RH$ is left morphic (as opposed to simply $x$ being left morphic in $RH$);   
nor do we know whether each finite subgroup $K$ of $G$ is contained in a finite subgroup $H$ of $G$ for which $RH$ is left morphic.
Certainly, both of these statements are trivially true if the group $G$ is a finite group (take $H = G$).      

We conclude with a few more examples.  
\begin{examples}  Let $R$ be a artinian principal ideal ring, and let $G$ be an infinite locally finite group for which $RH$ is a principal ideal ring for each nontrivial (finite) subgroup of 
$G$.  For instance, we may take $R$ to be a division ring, and $p$ a prime number which is invertible in $R$, we may take $G = \{x \in \CC: x^{p^r} = 1 \text{ for some } r \ge 0\}$.  Then, $RG$ is not a principal ideal ring 
by Passman's Theorem (it is infinite, but has no elements of infinite order), however, by \cite[Theorem 2.4]{morphgrp}, $RG$ is a left and right morphic ring. 
\end{examples}

\begin{examples} 
Let $R = \CC[t;\sigma]/(t^2)$, where $\sigma$ is an automorphism of $\CC$ which fixes the algebraic numbers
(there are such maps which are nontrivial, see, for instance, \cite{complex}).  
Then, any finite group $G$ is $R$-admissible (see the discussion preceding Theorem \ref{mainthm}).  In particular, 
$RG$ is a principal ideal ring for each finite group $G$.  
\end{examples}

\section{Appendix}

In this appendix, we will detail slight changes to the arguments found in \cite[Section 4]{passman} which allow one to replace the hypothesis that $K$ is a field with the hypothesis
that $K$ is a division ring in each of the results found in \cite[Section 4]{passman}.  

First, suppose $R$ is a ring which has a subring, with the same unity as $R$, which 
is isomorphic to $\Q$ or to $\Z_p$ for some prime $p$; we then view $R$ has having characteristic $0$ or characteristic $p$, accordingly.  
We may define the dimension subgroups as in \cite[Section 3.3]{passmanbook}, and 
\cite[Lemma 3.3.1]{passmanbook} and \cite[Lemma 3.3.2]{passmanbook} remain valid (we need $(x-1)^p = x^p - 1$, and we need to divide by positive integers, in the 
characteristic $p$ and $0$ cases, respectively).  Also, observe that $D_n(\M_n(R)) = D_n(R)$, since $\Delta(\M_n(R))^i = \M_n(\Delta(R)^i)$.      

We will next detail why \cite[Lemma 4.3]{passman} remains valid for division rings.  
In the proof of \cite[Lemma 4.3]{passman}, the first paragraph remains valid for any division ring $K$, with dimension interpreted as left $K$-vector space dimension.  The 
next paragraph (finding a subgroup $H$ for which $|H| \ne 0$ in the division ring) requires no changes.  By the properties cited for the dimension subgroups in this context, 
$G/H$ is a $p$-group, and if it is not cyclic, it has a homomorphic image which is elementary abelian of order $p^2$.  We need to make a slight change in the last paragraph, since the 
ring $KW$ need not be commutative.  However, arguing as in the first paragraph of the proof, if $\Delta(KW)$ is principal as a right ideal, say $\Delta(KW)= \alpha KW$, then 
$\Delta(KW) = KW \alpha$.  We conclude that $\Delta(KW)^p = KW \alpha^p$.  But, if $\alpha = \sum_{g \in W} a_g g$, then $\alpha^p = \sum_{g \in W} a_g^p = \epsilon(\alpha)^p = 0$.  
We conclude that if $\Delta(KW)$ is principal, it must be nilpotent of degree $p$.  It is, however, easy to see that the subring $\Z_p W$ of $KW$ is nilpotent of degree $2p-1>p$, from 
which it follows that $\Delta(KW)$ is not principal, and hence $G/H$ is a cyclic $p$-group.    

Finally, we observe that \cite[Lemma 4.4]{passman} remains valid for $K$ a division ring, with no changes needed.  
The proof of Passman's Theorem proceeds as before for the implications $(a) \Longrightarrow (b)$ and $(b) \Longrightarrow (c)$, using \cite[Lemma 4.3]{passman}
and \cite[Lemma 4.4]{passman} which hold for division rings.  The implication $(c) \Longrightarrow (a)$ is essentially unchanged, using Maschke's Theorem and \cite[Lemma 6]{sehgalfisher}
(more details of this type of argument are found in the proof of the main theorem in \cite{sehgalfisher}).

\section{Acknowledgements}  
I thank Don Passman for his helpful suggestions and correspondence regarding this material.  
I thank Alex Diesl for introducing me to the class of morphic rings, and also for 
many pleasant discussions regarding this material.  I also thank him for his many useful suggestions, in particular, for pointing out a simplification in the 
proof of Proposition \ref{prop3}.  It was also his suggestion to look at the ring $(\CC[t;\sigma]/(t^2)) C_3$ (see Example \ref{alex}), where $\sigma$ is complex conjugation (which together we showed was not a principal ideal ring directly) which helped eventually lead me to the main theorems in the case of associated graded rings.

\bibliography{PrincIdealGpRing}

\end{document}